\documentclass[12pt]{article}  
\usepackage{amsmath}
\usepackage{hyperref}
\usepackage{amsfonts}
\usepackage{amssymb}
\usepackage{graphics}
\usepackage{color}
\usepackage[pdftex]{graphicx}

\def\bee{\begin{equation}}
\def\eee{\end{equation}}

\def\FI{\scriptstyle{\rm FI}}
\def\hksqrt{\mathpalette\DHLhksqrt}
\def\DHLhksqrt#1#2{\setbox0=\hbox{$#1\sqrt{#2\,}$}\dimen0=\ht0
\advance\dimen0-0.2\ht0
\setbox2=\hbox{\vrule height\ht0 depth -\dimen0}%
{\box0\lower0.4pt\box2}}
\begin{document}

\thispagestyle{empty}
\bigskip
\centerline{    }
\centerline{\Large\bf Continued fractions constructed}
\centerline{\Large\bf from prime numbers}
\bigskip\bigskip
\centerline{\large\sl Marek Wolf}

\begin{center}
e-mail: primes7@o2.pl
\end{center}

\bigskip\bigskip

\begin{center}
{\bf Abstract}\\
\end{center}

\begin{minipage}{12.8cm}
We give 50 digits values of the simple continued fractions  whose denominators
are formed from  a) prime numbers, b) twin primes, c) generalized $d$-twins,
d) primes of the form $m^2+n^4$, e)primes of the form $m^2+1$,
f) Mersenne primes and g) primorial primes. All these continued fractions belong to the set of
measure zero of exceptions to the theorems of Khinchin and Levy. We claim
that all these continued fractions are transcendental numbers.
Next we propose the conjecture which indicates the way to deduce the transcendence
of some continued fractions from transcendence of another ones.
\end{minipage}

\bigskip\bigskip

\section{Introduction}

Let  $a_0$ be an integer and let  $a_k,~k=1, 2, \ldots, n$ are
positive integers (in general $a_k$ can be arbitrary complex numbers, see e.g.
\cite{Khovanskii}). Then
\bee
r=[a_0; a_1, a_2, a_3, \ldots, a_n]\equiv
a_0+\cfrac{1}{a_1 + \cfrac{1}{a_2 + \cfrac{1}{a_3+\ddots\cfrac{1}{a_n}}}}
\eee
is the simple (i.e. with all nominators equal to 1)
finite continued fraction. The numbers $a_k,~k=1, 2, \ldots, n$ are called
partial quotients and
\bee
\frac{P_k}{Q_k}=[a_0; a_1, a_2, a_3, \ldots, a_k], ~~~~~~~~k=1, 2, \ldots n
\eee
is called  the $k$-th convergent of $r$. If for the infinite continued fraction
\bee
[a_0; a_1, a_2, a_3, \ldots]
\label{ulamek}
\eee
the sequence of  convergents  $P_n/Q_n$ converges to some
limit $r$ when $n\rightarrow \infty$ then we say that the infinite continued fraction
$[a_0; a_1, a_2, a_3, \ldots]$ is equal to $r$.
The  convergence of the continued fraction (\ref{ulamek}) is linked to  the
behavior   of the sum of partial quotients  $a_n$:
\bee
{\rm sequence~~}\frac{P_n}{Q_n} {\rm ~~is ~convergent~to}~~r
~~\Leftrightarrow ~~\sum_{n=1}^\infty a_n {\rm~~ is ~ divergent}
\label{theorem}
\eee
see e.g. \cite[Theorem 10, p.10]{Khinchin}. It means that for convergence
of the continued fraction it is necessary that both $P_n, Q_n \to \infty$ in such a way,
that the ratio $P_n/Q_n$ has a definite limit for $n\to \infty$.
If  the infinite continued fraction is convergent then the values of the
convergents  ${P_k(r)/Q_k(r)}$ approximate the value of $r$  with accuracy
at least $1/Q_k Q_{k+1}$
\cite[Theorem 9, p.9]{Khinchin}:
\bee
\left|r - \frac{P_k}{Q_k} \right |<\frac{1}{Q_k Q_{k+1}}<\frac{1}{Q_k^2a_{k+1}}<\frac{1}{Q_k^2}.
\label{error}
\eee
Rational numbers have finite continued fractions, quadratic irrationals
have periodic infinite continued fractions and vice versa: eventually periodic
continued fractions represent quadratic  surds. All remaining irrational
numbers have non-periodic continued fractions.

Khinchin has proved that \cite[p.93]{Khinchin}
\bee
\lim_{n\rightarrow \infty} \big(a_1 a_2 \ldots a_n\big)^{\frac{1}{n}}=
\prod_{m=1}^\infty {\left\{ 1+\frac{1}{m(m+2)}\right\}}^{\log_2 m} \equiv K_0 \approx 2.685452001\ldots
\label{Khinchin}
\eee
is a constant for almost all real $r$, see also \cite{Ryll-Nardzewski1951},
\cite[\S 1.8]{Finch}.  The exceptions are  rational numbers,
quadratic irrationals and some irrational numbers too, like for example the
Euler constant $e=2.7182818285\ldots$, but this set of exceptions is of the
Lebesgue measure zero. The constant $K_0$ is called the Khinchin constant.

In 1935 Khinchin \cite{Khinchin-1935} has proved that for almost all real   $r$ the
denominators $Q_n(r)$ of the convergents of the  continued fraction expansions for $r$
satisfy  $\lim_{n\to \infty} \sqrt[n]{Q_n(r)}=L_0$ and in  1936 Paul Levy
\cite{Levy1936} found an explicit expression for this constant $L_0$:
\bee
\lim_{n \to \infty}\sqrt[n]{Q_n(r)}=e^{\pi^2/12\log(2)}\equiv L_0=3.27582291872\ldots
\label{Levy}
\eee
All presented below continued fractions belong to this exceptional set of
irrationals for which the geometric means of the denominators
$(a_1 a_2 \ldots a_n)^{1/n}$  and the $n$-th radical roots of the denominator
$Q_n^{1/n}$  tend to infinity.

In this paper we will consider continued fractions with partial quotients given
by  an infinity of all primes as well as primes of special form belonging to families
containing conjectured infinity of members. All these  continued fractions
are non-periodic, and thus are
irrational, but we claim that all of them are also {\it transcendental}.
In Sect. 3 we review some facts and theorems concerning the transcendentality
of continued fractions. In Sect. 4 some experimental results regarding
transcendentality of numbers constructed from primes are presented.

\section{Seven examples}

In consecutive sections we will discuss  the following cases:
the set of all primes $2, 3, 5, 7, \ldots$,  twin primes, generalized $d$-twins, i.e.
pairs of adjacent primes separated by $d$,  primes of the form  form $m^2+n^4$,
primes given by the quadratic form $m^2+1$. Next are considered sparse Mersenne primes
and at the end even sparser primorials  primes.

It is possible to consider other families of primes, like Sophie Germain primes
(it is conjectured that there are infinitely many of them),
irregular primes (Jensen in 1915 proved that there are infinitely many
of them), regular primes of which it was  conjectured  that $e^{-1/2}\approx 61\%$
of all prime numbers are regular, the Cullen numbers $n2^n+1$ when they are primes etc.  but we leave it for further studies.

Except  Sections 2.1  and 2.4, where we will  treat the  families of primes
containing rigorously proved  infinity of members, all remaining consideration
are performed under the assumption  there is infinity  of  primes in each class
of primes, although proofs of infinitude of all these sets of primes seems
to be very far away. Thus many of our reasonings are heuristical.

The examples are in order of sparseness of each family of primes.

\subsection{The set of all primes}

Let us put $a_n=p_n$ where $p_n$ denotes the $n$-th primes:
$[0; 2, 3, 5, 7, 11, 13, \ldots]$. As there is
an infinity of primes the condition (\ref{theorem}) is fulfilled and let
us denote the limit of the continued fraction by
\bee
u=[0;2, 3, 5, 7, 11, 13, \ldots] = \cfrac{1}{2+\cfrac{1}{3+\cfrac{1}{5+\cfrac{1}{7+\cfrac{1}{11+\ddots}}}}}
\label{def-up}
\eee
Using PARI system \cite{PARI} and all 1229 primes up to 10000 it is possible to obtain
over 8000 digits of the above continued fraction in just a few seconds because
\bee
[0;2, 3, 5, 7, 11, 13, \ldots, 9973]=\frac{3.38592889\ldots \times 10^{4297}}{7.83177791\ldots \times 10^{4297}}
\eee
and the product of $Q_k Q_{k+1}$ on the rhs of (\ref{error}) is larger than $10^{8500}$.
The  first 50 digits of $u$  reads:
\bee
u=0.43233208718590286890925379324199996370511089688\ldots .
\label{u-1}
\eee
This  number is not recognized at the  Symbolic Inverse Calculator
(http://pi.lacim.uqam.ca/eng/) maintained by Simone Plouffe. Accidentaly, it
is very close to the one of Renyi's parking constants $m_{\rm R}=(1-e^{-2})/2=0.43233235838\ldots$,
see  \cite[pp. 278--283]{Finch}: $m{\rm R}-u=2.712\ldots\times10^{-7}$.

It is possible to obtain analytically the geometrical means of the partial quotients
in (\ref{def-up}). It is well known (see e.g. \cite[Chap.4] {Ellison}), that
the Chebyshev function $\theta(x)$  behaves like:
\bee
\theta(x) \equiv  \sum_{p\leq x} \log(p) = x + \mathcal{O}(\sqrt{x}).
\label{Chebyshev}
\eee
Thus skipping the error term we have
\bee
\prod_{k=1}^n p_k = e^{p_n}.
\eee
It is well known that \cite[Sect. 2.II.A]{Ribenboim} that
\bee
p_n=n \log(n) + n(\log \log (n)-1) + o\left(\frac{n\log \log (n)}{n}\right).
\eee
For our purposes it suffices to know that
\bee
p_n > n \log(n) ~~~~ {\rm for ~~} n>1
\eee
see e.g. \cite{Schoenfeld1962}. Hence we can write
for the geometrical means of the partial quotients the estimation:
\bee
\big(a_1 a_2 \ldots a_n\big)^{\frac{1}{n}}= \left(\prod_{k=1}^n p_k \right)^{\frac{1}{n}}
= \left(e^{p_n}\right)^{\frac{1}{n}} > n \rightarrow \infty
\label{u-GMK}
\eee
thus the continued fraction $u$ belongs to the set of measure zero of exceptions
to the Khinchin Theorem (\ref{Khinchin}).  It is also an exception to the Levy Theorem,
because from the general properties of continued fractions:
\bee
Q_{n+1} =a_nQ_{n}+Q_{n-1}
\eee
we have $Q_n>\prod_{k=1}^n p_k > n^n$ and thus $Q_n^{1/n}  \to \infty$ in
contrast to (\ref{Levy}). It is an explicit example of the continued fraction
with unbounded $(Q_n)^{1/n}$.

\subsection{Twin primes}

The twin prime conjecture states that there are infinitely many pairs of primes
$(t_n, t_{n+1})$ differing by two: $t_{n+1}-t_n=2$. Let $\pi_2(x)$ denote the number
of pairs of twin primes  $(t_n, t_{n+1}) $ smaller than $x$. Then
the conjecture B of Hardy and Littlewood \cite{Hardy_and_Littlewood} on
the number of prime pairs $p, p+d$ applied to the case $d=2$ gives, that
\bee
\pi_2(x) \sim  C_2 \int_2^x \frac{u}{\log^2(u)} du = C_2 \frac{x}{\log^2(x)}+\ldots,
\label{conj}
\eee
where $C_2$ is called ``twin constant'' and is defined by the following
infinite product:
\begin{equation}
C_2 \equiv 2 \prod_{p > 2} \biggl( 1 - \frac{1}{(p - 1)^2}\biggr) =
1.32032363169\ldots
\label{stalac2}
\end{equation}
If there is indeed (as everybody believes, see e.g. \cite{Koreevar}) an infinity
of twins, then
the continued fraction
\bee
u_2=[0; 3, 5, 5, 7, 11, 13, 17, 19, \ldots]
\label{def-u2}
\eee
should be   infinite, non-periodic and  convergent.
We count here 5 two times as it is a customary  way
of defining the Brun's constant \cite{Wrench1974} and it an only
case of double appearance of a prime in the set of twins  as for adjacent twin
pairs $(p-2, p)$ and $(p, p+2)$  one of numbers $(p-2, p, p+2)$ always
is divisible by 3.
Again performing calculations in PARI and using primes $<10000$ we found here 205
twin pairs (but only 409 different primes) and first 50 digits of the continued
fraction (\ref{def-u2}) are
\bee
u_2= 0.31323308098694591263078648647217280043925117451\ldots .
\eee
There is much less terms in $u_2$ up to $10000$ than  primes $<10000$ in $u$,
hence the value of
$u_2$ was obtained with accuracy about 2900 digits.  We have checked using Plouffe's
Symbolic Inverse Calculator (http://pi.lacim.uqam.ca/eng/), that this  constant
is not recognized as a combination of other mathematical quantities.

Because twin primes are sparser than all primes we have $t_n>p_n$ thus in view
of (\ref{u-GMK}) the geometrical means  $(3\cdot 5 \ldots t_n)^{1/n}$ will diverge
even faster, hence the  continued fraction $u_2$ belongs to the set of exceptions to the
Khinchin Theorem. It is also a counterexample to the Levy Theorem, because
denominator $Q_n(u_2)$ of the $n$-th convergent of  $u_2$ is larger than the
denominator  $Q_n(u)$ of the  $n$-th convergent of $u$.

\subsection{Generalized $d$-twins}

It is natural to consider the whole family of continued fractions
$u_d$,  $d=2, 4, 6, 8, \ldots$ formed from the {\it consecutive} primes separated by $d$:
$p_{n+1}-p_n=d$. We put this example here after twins, although for sufficiently
large  $d$ the  primes  $p_{n+1}-p_n=d$  will be even sparser than say Mersenne primes
and from the other side $d=2$ are less frequent  than $d=6$, see \cite{champions}.
The consecutive primes separated by $d=4$ are sometimes called Cousins,
\cite{Wolf-RW}.  For example, in the case of $d=6$ we have :
\[
u_6=[0; 23, 29, 31, 37, 47, 53, 53, 59, 61, 67, 73, 79, 83, 89, 131, 137, 151, 157, 157, 163, ...]
\]
and some primes $p_n$ when $p_n-p_{n-1}=p_{n+1}-p_n$  do appear twice (in the case of
$u_2$ only 5 appears two times).
As in the case of twins it is conjectured that for each $d$  there is an infinity
of prime pairs $(p_{n+1}, p_n)$   with  $p_{n+1}-p_n=d$,
see e.g. \cite{Brent1974}, \cite{champions}. From this conjecture it follows that
 the numbers $u_d$ are irrational. Using PARI/GP we have calculated
the values of $u_d$ up to $d=570$, what took four days  of CPU time on the  64 bits
AMD Opteron 2700 MHz processor. We have searched for primes up to
$2^{44}=1.759\ldots \times 10^{13}$ and the largest encountered gap between
consecutive primes was $d=706$ which appeared only once. We have
calculated $u_d$  if there was a number of  gaps of given $d$ sufficient to
determine  $u_d$ with at least a few hundreds digits (for example, there were
17 pairs of 570--twins up to $2^{44}$). The Table I gives a sample
of obtained values with 50 digits accuracy. The whole file with
275 values of $u_d$ given with more than 110 digits is available from
the author webpage  http://www.ift.uni.wroc.pl/$\sim$mwolf/u\_d.dat.

For large $d$ the value of $u_d$ is practically determined by the first occurrence
$p_f(d)$ of that gap ---  pairs of  consecutive primes
with gap  $d\gg 2 $ are separated by very large intervals, for example first $d=540$
appears between $(738832927927,   738832928467)$ and next gap $d=540$ is
between $(3674657545087,  3674657545627)$. It was conjectured by Shanks in 1964
\cite{Shanks1964} that the
gap $d$ appears for the first time at the prime $p_f(d)\sim e^{\sqrt{d}}$. We have
given heuristic arguments  \cite{Wolf-first-d}   that
\bee
p_f(d)\sim \sqrt {d} \exp\left({1\over 2}\sqrt{\ln^2(d) + 4 d}\right)
\eee
and for $d\gg 1$ simply
$p_f(d)\sim \sqrt{d} e^{\sqrt{d}}$. Thus we claim that for large $d$
there should be the approximate formula:
\bee
u_d \approx [0;  \sqrt{d}e^{\sqrt{d}},  \sqrt{d}e^{\sqrt{d}}+d]\approx \frac{1}{ \sqrt{d}e^{\sqrt{d}}}.
\eee

The plot of $u_d$ and comparison with the Shanks and our conjecture is
given in the Fig.1.

Again like $u$ and $u_2$ continued fractions $u_d$ belongs to the set of
exceptions to the Khinchin Theorem and Levy Theorem.

\subsection{Primes of the form $m^2+n^4$}

In the seminal paper \cite{Iwaniec-Friedlander} John Friedlander and  Henryk Iwaniec
have proved that there exists infinity of primes of the form $m^2+n^4$.
More precisely, if $\pi_{\FI}(x)$ denotes the number of primes of the form
$m^2+n^4<x$ then approximately
\bee
\pi_{\FI}(x) \sim \frac{C_{\scriptstyle{ \rm FI}}x^{3/4}}{\log(x)}
\label{Iwaniec}
\eee
where the constant $C_{\FI}=\hksqrt{2}\Gamma(\frac{1}{4})^2/3\pi^{3/2}=
1.112835788898764\ldots$ and here $\Gamma$ is the Euler Gamma function.
Thus taking as partial quotients of the continued fraction primes of the form
$m^2+n^4$ for sure we will obtain an irrational number  which we will
denote $u_{\rm \FI}$:
\bee
u_{\rm FI}=[0; 2, 5, 17, 17,  37, 41,  97, 97, \ldots]
\eee
Like in previous examples some primes appear twice: $17=4^2+1^4=1^2+2^4,
~ 97=9^2+2^4= 4^2+3^4$ etc. Looking for  all primes of this form with
$1\leq m \leq 100$ and $1\leq n \leq 10$ (the largest prime was $19801=99^2+10^4$)
we get the value  of $u_{\FI}$ with over 1100
digits accuracy; the first 50 digits of it are:
\bee
u_{\rm FI}=0.455024816490170022369052808279744824105755548905\ldots
\eee
Let us notice that
\[
1/(2 + 1/(5 + 1/(17 + 1/(17 + 1/(37 + 1/(41 + 1/(97 + 1/98)))))))=
\]
\[
\frac{20993638525}{46137348479}=0.455024816490170022369048157801049432084768331968\ldots
\]
and the difference between this value and $u_{\rm FI}$ is less than $10^{-23}$ !

\begin{figure}
\vspace{-2.3cm}
\hspace{-3.5cm}
\begin{center}
\includegraphics[width=\textwidth,angle=0]{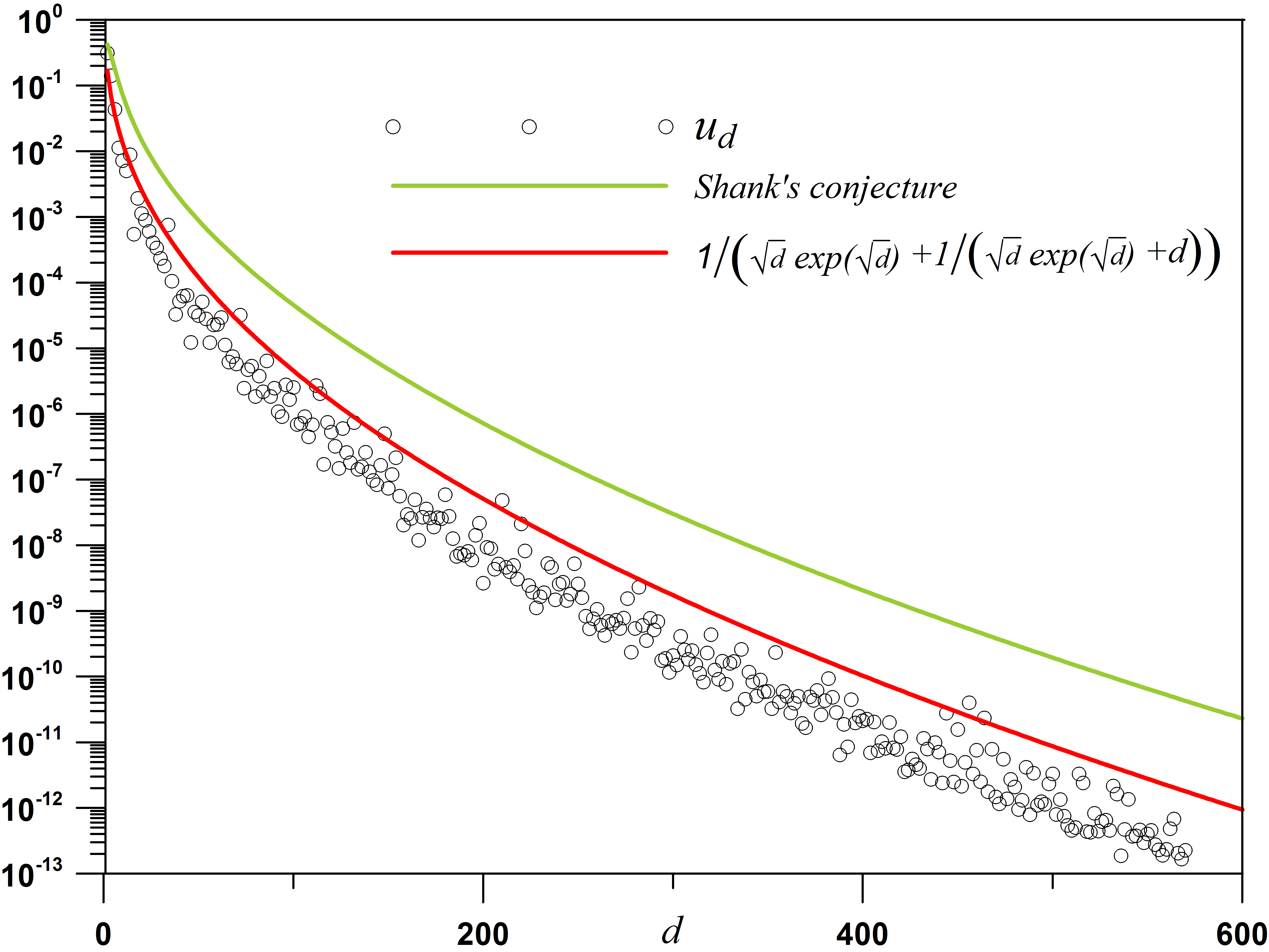} \\
\vspace{0.2cm}{\small  Fig.1  The plot of $u_d$ and two approximations:  in green
the Shank's conjecture $\cfrac{1}{e^{\sqrt{d}}+\cfrac{1}{e^{\sqrt{d}}+d}}
$ and in red
our conjecture
$\cfrac{1}{\sqrt{d}e^{\sqrt{d}}+\cfrac{1}{\sqrt{d}e^{\sqrt{d}}+d}}$.}
\end{center}
\end{figure}

\newpage
\noindent
\begin{center}
~~~~~~~~~~~~~~~~~~\\
\vspace{-2cm}
Table I \\
\medskip
\begin{tabular}{|c|c|} \hline
$  d  $ & $ u_d $   \\  \hline
4 &     1.4103814184127409729946079947661391024642878552250$\times 10^{-1}$ \\ \hline
6 &     4.3413245800886640441937906138426444157119875018764$\times 10^{-2}$ \\ \hline
8 &     1.1234653732060451418609230935360294984983811524705$\times 10^{-2}$ \\ \hline
10 &     7.1938972705064358418419102215951120335820544247877$\times 10^{-3}$ \\ \hline
12 &     5.0250059564863844924667112008186998625931272954692$\times 10^{-3}$ \\ \hline
14 &     8.8489409307271044901495673780577102976304420791245$\times 10^{-3}$ \\ \hline
16 &     5.4614948350881467294308534284337241698766002935218$\times 10^{-4}$ \\ \hline
18 &     1.9120391314299159400657740968697274305281924125799$\times 10^{-3}$ \\ \hline
20 &     1.1273943145526585257207207582991176443515379616999$\times 10^{-3}$ \\ \hline
22 &     8.8573891094929851372874303656530678911673854053699$\times 10^{-4}$ \\ \hline
24 &     5.9916096230989554005997263265407846890656053212565$\times 10^{-4}$ \\ \hline
26 &     4.0371410525148524468010569219212401713453876041188$\times 10^{-4}$ \\ \hline
28 &     3.3658696996531260967017397551173798914121748535404$\times 10^{-4}$ \\ \hline
30 &     2.3272049015980164345521674554989676374011829679698$\times 10^{-4}$ \\ \hline
32 &     1.7885887465418665415382499015390795012182483537844$\times 10^{-4}$ \\ \hline
34 &     7.5357908538425634007656299916322144807040843935028$\times 10^{-4}$ \\ \hline
36 &     1.0470107727765143055064789951193220804598138780293$\times 10^{-4}$ \\ \hline
38 &     3.2687215994929278130910770751451289367042590019431$\times 10^{-5}$ \\ \hline
40 &     5.1725029603788623137563671868924142637218125718293$\times 10^{-5}$ \\ \hline
42 &     6.1954029872477528100134249220879079074519481392595$\times 10^{-5}$ \\ \hline
44 &     6.3763310332564890009355447509046689625278954819441$\times 10^{-5}$ \\ \hline
46 &     1.2275511580096446939755547564625149207372813752259$\times 10^{-5}$ \\ \hline
48 &     3.5424563347877377649245656903453981296411399487963$\times 10^{-5}$ \\ \hline
50 &     3.1341084997626641267187094247975857118584579732840$\times 10^{-5}$ \\ \hline
\vdots & \vdots \\ \hline
566 &     2.0417988154535953561248601983565125430801657124094$\times 10^{-13}$ \\ \hline
568 &     1.6638019955234637865242752590874891539355008513604$\times 10^{-13}$ \\ \hline
570 &     2.2511824714719308536000694530283450847909292429681$\times 10^{-13}$ \\ \hline
\end{tabular}
\end{center}

\subsection{Primes of the form $m^2+1$}

Now let us consider the set of prime numbers
\bee
\mathcal{Q}=\{2,5,17,37,101,197, 257, 401, 577, 677, 1297, 1601,   \ldots\}
\label{set_Q}
\eee
given by the quadratic polynomial  $m^2+1$ and let $q_n$ denote the $n$-th prime
of this form.  By  the conjecture E of Hardy and Littlewood \cite{Hardy_and_Littlewood}
the number $\pi_q(x)$  of primes  $ q_n< x$ of the form
$q_n = m^2 + 1$ is  given by
\bee
\pi_q(x) \sim C_q  \frac {\sqrt{x}}{\log(x)},
\label{conj-E}
\eee
where
\bee
C_q=\prod_{p\geq 3}\biggl(1-\frac{(-1)^{(p-1)/2}}{p-1}\biggr) =
1.372813462818246009112192696727\ldots
\eee
Comparing it with (\ref{Iwaniec}) we see that indeed primes $m^2+1$ are sparser
than primes $m^2+n^4$. For example up to $10^8$ there are 65162 primes of the
form $m^2+n^4$ and only 841 primes of the form $m^2+1$.
Although  the conjecture (\ref{conj-E}) remains unproved there is no doubt in its validity.
Thus let us create the presumedly infinite  continued fraction by identifying
$a_n=q_n, n\geq 1$:
\bee
u_q=[0; 2,5,17,37,101,197, 257, 401, 577,  677, 1297, 1601, \ldots] .
\label{def-uq}
\eee
Using 841 primes of the form $m^2+1$ smaller than $10^8$ and performing the calculations
in PARI with precision set to 20000 digits
we get over 11000 digits of $u_q$ as the ratio on the rhs of (\ref{error})
was $<10^{-11700}$. First 50 digits of $u_q$ reads:
\bee
u_q=0.45502569980199468718020210263808421898137687948\ldots
\eee
Let us remark that $u_{\FI}-u_q=8.833\ldots \times 10^{-7}$.

There is no known formula analogous to (\ref{Chebyshev}) for primes of the form
$m^2+1$, but because $q_n\geq p_n$ the geometrical means of $2\cdot 5 \cdot 17 \ldots q_n$
will diverge faster than (\ref{u-GMK}). It is possible to obtain very rough speed
of divergence of $(2\cdot 5 \cdot 17 \ldots q_n)^{1/n}$. Namely,
making use of (\ref{conj-E})  and   inverting $\pi_q(q_n)=n$ we get:
\bee
q_n\sim \left(\frac{2 n \log(n/C_q)}{C_q}\right)^2 +
2\log\left(\frac{n}{C_q}\right)\log\log\left(\frac{n}{C_q}\right)
\eee
Because $2>C_q$ it follows that $2\cdot 5 \cdot 17 \ldots q_n$ grows faster
than $2^{2n}(n!)^2/C_q^{2n}>(n!)^2$
and the Stirling formula for $n!$ gives that $(2\cdot 5 \cdot 17 \ldots q_n)^{1/n}$
grows faster than $n^2$ and again $u_q$ is the exception to the Khinchin Theorem
as well as to the Levy Theorem.

\subsection{Mersenne primes}

The Mersenne primes $\mathcal{M}_n$ are the  primes of the form $2^p-1$ where $p$
must be a prime, see  e.g. \cite[Sect. 2.VII]{Ribenboim}. Only 47 primes of this form
are currently known, see Great Internet Mersenne Prime Search (GIMPS)
at www.mersenne.org.    For many years the largest known primes
are the Mersenne primes, as the Lucas--Lehmer primality test (applicable only to
$\mathcal{M}_n = 2^p-1$) needs just a multiple of $p$ steps, thus the complexity
of checking primality of $\mathcal{M}_n$ is $\mathcal{O}(\log(\mathcal{M}_n))$.
Let us remark that algorithm of Agrawal, Kayal and Saxena (AKS) for general
prime $p$  works in  about $\mathcal{O}(\log^{7.5}(p))$ steps and
modification by Lenstra and Pomerance in about $\mathcal{O}(\log^{6}(p))$  steps.

Again there is no proof of
the infinitude of $\mathcal{M}_n$ but a common belief is  that as there are
presumedly infinitely many even perfect numbers thus there is also an infinity
of Mersenne primes.

\begin{figure}
\vspace{-2.3cm}
\hspace{-3.5cm}
\begin{center}
\includegraphics[width=0.4\textheight,angle=0]{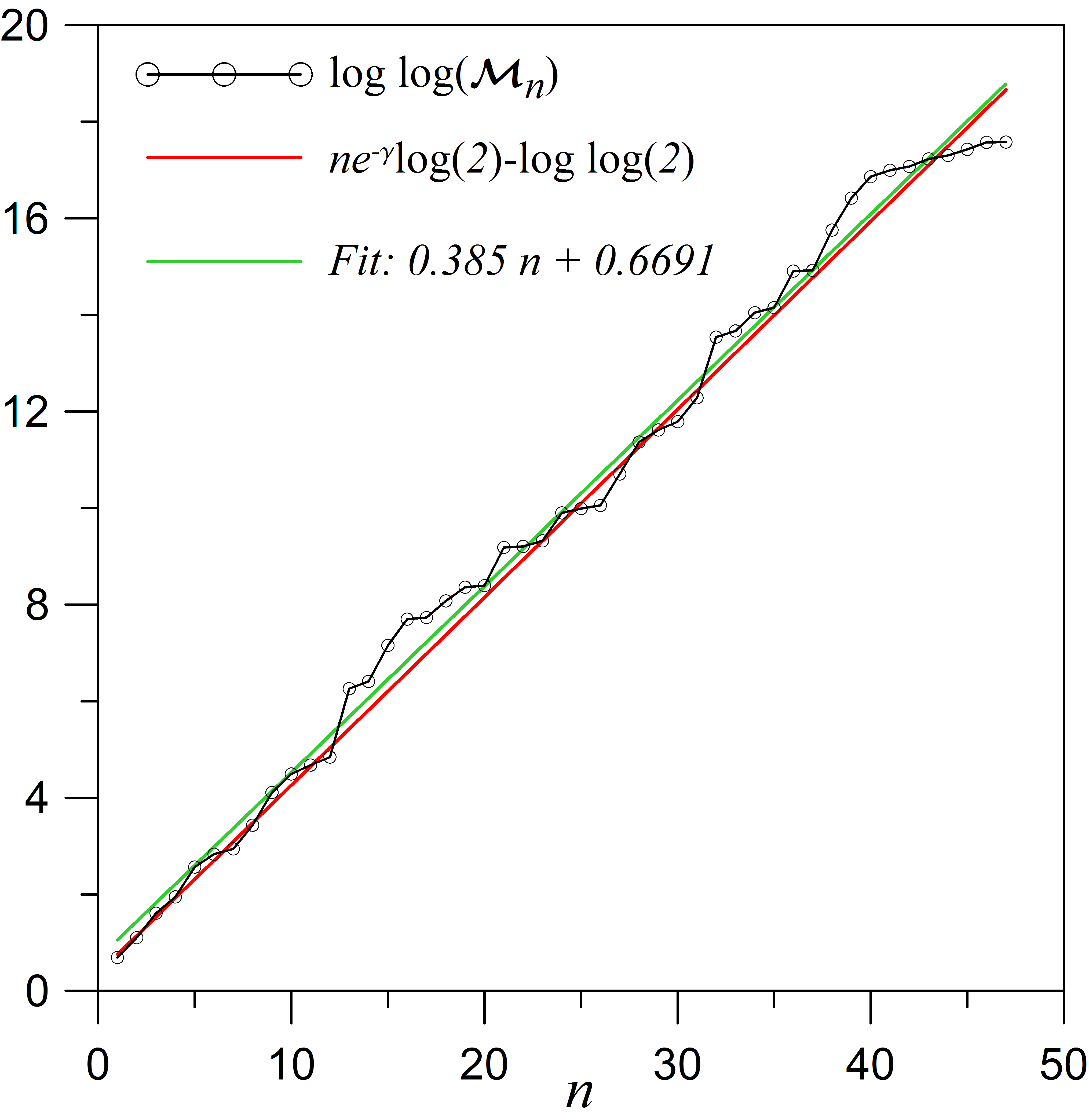} \\
\vspace{0.2cm} {\small Fig.2  The plot of $\log\log(\mathcal{M}_n)$ and the Wagstaff conjecture
(\ref{Wagstaff}).  The fit was made to all known $\mathcal{M}_n$ and it is
$0.3854n + 0.6691$, while $n e^{-\gamma}\log(2)n -\log\log(2) \approx
0.3892 n +0.3665$. The rather good  coincidence of $\log\log(\mathcal{M}_n)$
and (\ref{Wagstaff}) is seeming, as to get original $\mathcal{M}_n $'s the errors
are amplified to huge values by double exponentiation.}
\end{center}
\end{figure}

Let us define the supposedly infinite and convergent continued fraction
$u_{\mathcal{M}}$ by taking $a_n=\mathcal{M}_n$:
\bee
u_{\mathcal{M}}=[0; 3,~ 7,~ 31,~ 127,~ 8191,~ 131071,~ 524287,~ 2147483647,~ \ldots]
\label{def-uM}
\eee

Using all 47 Mersenne primes  $3, 7, 31, \ldots, 2^{43112609}-1$
in a couple of minutes we have calculated $u_{\mathcal{M}}$
with the precision better than  $10^{-121949117}$;  first
50 digits of $u_{\mathcal{M}}$ are:
\bee
u_{\mathcal{M}}=0.31824815840584486942596202748140694243806236564\ldots
\eee

Of course $u_{\mathcal{M}}$ is also the exception to the Khinchin and Levy Theorems
in view of the very fast growth of $u_{\mathcal{M}}$ --- Wagstaff conjectured \cite{Wagstaff1983},
that ${\mathcal{M}_n}$ grow  doubly exponentially:
\bee
\log_2 \log_2 \mathcal{M}_n \sim n e^{-\gamma}.
\label{Wagstaff}
\eee
where $\gamma=0.57721566\ldots$ is the Euler--Mascheroni constant.
In the Fig. 2 we compare the Wagstaff conjecture with all 47 presently known
Mersenne primes.

\subsection{Primorial primes}

If $p_n$ is the $n-$th prime number then numbers  of the form
$2\times 3\times 5 \cdots \times p_n\equiv p_n \sharp $
are  called primorials and $\sharp $ stands  here by analogy of exclamation mark in
the factorial. The primorials are expressed directly by the Chebyshev function
$\theta(x)$:
\bee
p_n\sharp = e^{\theta(p_n)}=e^{(1+o(1))p_n}.
\eee
For some primes $r_n$ the numbers $r_n\sharp \pm 1$ are primes. They are called
primorial primes and are even sparser than Mersenne primes as we will see below.
Despite this rareness of primorial primes it was conjectured that there is infinity
of them \cite{Caldwell-2002} and thus the continued fractions
\bee
u_{r+}=[0; 3 ,  7,  31, 211,  \ldots ]
\eee
obtained from primorial primes of the form $r_n\sharp + 1$ will be at least irrational
number, as   well  as the continued fraction  build from primorial primes of the
form $r_n\sharp - 1$:
\bee
u_{r-}=[0; 5, 29, 2309, \ldots ].
\eee
From the known presently only 22 (see sequence A005234 in OEIS) primorial primes
$r_n \sharp + 1$ we get the continued fraction
\bee
u_{r+}=[0; 3 ,  7,  31, 211,  \ldots,   42209\sharp +1, 145823\sharp+1,
366439\sharp+1 	, 392113\sharp+1]=
\eee
\[
0.318248165083690124777685589996787844788657122331533049467\ldots
\]
with the error less than $10^{-914474}$. Let us remark that $u_{r+}-u_{\mathcal{M}}=
6.678\ldots \times 10^{-9}$, although only three first primes (3, 7, 37)  are the same.

From all 18 presently known
(see sequence  A006794 in  OEIS)  primes of the form $r_n \sharp - 1$  we get
\bee
u_{r-}=[0; 5, 29, 2309, \ldots,   15877\sharp - 1]=
\eee
\[
0.198630157303503810875201233614346862875870630898479777625647\ldots
\]
with the error less than $10^{-48415}$.
Chris K. Caldwell  and Yves  Gallot gave heuristic arguments \cite{Caldwell-2002}
that there is infinity of primorial primes of both kinds. More precisely,
they claim that  the expected  numbers  of primorial primes of each of the forms
$r\sharp \pm 1$  with  $r < x$ are both approximately $e^\gamma \log(x)$.
From $n=e^\gamma \log(r_n)$ we get that
\bee
r_n \sim e^{n/e^\gamma},
\eee
where $r_n$ stands for $n$-th prime giving the primorial prime of the form
$r_n \sharp \pm 1$. Then the $n$-th primorial prime, and hence $a_n$ of $u_{r \pm}$,
will be
\bee
a_n=r_n \sharp =\exp\left(e^{-\gamma}\sum_{p\leq r_n} p\right).
\eee
From the formula
\bee
\sum_{p\leq n} f(p)=\int_2^n \frac{f(x)}{\log(x)} dx + f(2)li(2)+f(n)\left(\pi(n)-li(n)\right)
-\int_2^n \{\pi(x)-li(x)\}dx
\label{suma-po-p}
\eee
where $\pi(x)$ is the number of primes $<x$ and $li(x)=\int_2^x dx/\log(x)$ is the
logarithmic integral,  see  \cite[eq.(2.26)]{Schoenfeld1962} we get
\bee
\sum_{p\leq n} p = li(n^2)+ error= \frac{n^2}{2\log(n)}+error'.
\eee
Here $error'$ besides expressions on rhs in (\ref{suma-po-p}) contains also higher
terms coming from the asymptotic expansion of $li(x^2)$.  Finally we obtain
\bee
r_n\sharp \approx \exp\left(\frac{1}{2n}\exp\left(2e^{-\gamma}n\right)\right).
\eee
From this it follows that the $n$-th primorial prime is much larger than
the $n$-th Mersenne prime  $\mathcal{M}_n\sim 2^{2^{n/e^\gamma}}$. Indeed, the ratio:
\bee
\frac{\log(r_n\sharp)}{\log(\mathcal{M}_n)}=\frac{1}{2\log(2) n}\left(\frac{e^2}{2}\right)^{ne^{-\gamma}}
\eee
grows with $n$.

\section{Continued fractions and transcendence}

There is a vast literature concerning the transcendentality of continued fractions.
The Theorem of H. Davenport and K.F. Roth \cite{Davenport1955} asserts, that if
the denominators $Q_n$ of convergents of the continued fraction
$r=[a_0; a_1, a_2, \ldots]$ fulfill
\bee
\limsup_n \frac{\sqrt{\log(n)} \log (\log (Q_n(r)))}{n}=\infty
\label{Davenport-Roth}
\eee
then $r$ is transcendental. This  theorem requires for the transcendence
of $r$ very fast increase of denominators of the convergents: at least doubly
exponential growth
is required for \ref{Davenport-Roth}. The set of continued fractions which can satisfy
the  Theorem of H. Davenport and K.F. Roth is of measure zero, as it follows from the
Theorem 31 from the Khinchin's book \cite{Khinchin}, which asserts there exists an absolute
constant $B$ such that for {\it almost all} real numbers $r$ and sufficiently large
$n$ the denominators of its  continued fractions  satisfy:
\bee
Q_n(r)<e^{Bn}.
\eee

The paper of A. Baker \cite{Baker1962} from 1962 contains a few theorems
on the  transcendentality of  Maillet type  continued fractions \cite{Maillet},
i. e.  continued fractions with bounded partial quotients which have transcencendental
values. In the paper \cite{Adamczewski2007} B.  Adamczewski and  Y. Bugeaud,
among others, have improved (\ref{Davenport-Roth}) to the form:
\bee
{~~~~\rm If}~~~\limsup_n \frac{ \log (\log (Q_n(r)))}{n^{2/3}(\log(Q_n(r)))^{2/3}\log(\log(Q_n(r)))}=\infty
\label{Adamczewski-Bugeaud}
\eee
then $r$ is transcendental.

Besides Maillet  continued fractions there are some specific families of
other continued fractions of which it is known that they are
transcendental. In the papers \cite{Queffelec}, \cite{Adamczewski2007-AMM} it was
proved that the Thue--Morse
continued fractions with bounded partial  quotients are transcendental.
Quite recently there appeared the preprint
\cite{Bugeaud2010} where the transcendence of the Rosen continued fractions was
established. For more examples see \cite{Adamczewski-Bugeaud-Davison}.

Taking as the partial quotients $a_n$ different sequences of numbers
leads to real numbers which  very often turn out be transcendental.
For example, the continued fraction $s$ for which $a_n=n$:
\bee
s =[0; 1, 2, 3, 4, \ldots]=1/(1+1/(2+1/(3+1/(4+\ldots)))) = 0.697774657964\ldots
\label{natural}
\eee
is transcendental. Let us mention that the continued fraction $s'$ with all
partial quotients equal to consecutive odd numbers:
\bee
s'=[0; 1, 3, 5, 7, 9, \ldots]=-i\tan(i)=\frac{e-e^{-1}}{e+e^{-1}}=
0.761594155955764888119\ldots
\label{odd}
\eee
is also transcendental. All these facts are special cases of the results obtained by
Carl Ludwig Siegel in 1929 in a long paper  \cite{Siegel-1929}.
In particular, the continued fractions $s_D=[0; a_1, a_2, a_3, \ldots]$
are transcendental when the $a_n$'s are rational and form an arithmetical sequence
of the difference $D$ and first element $A: ~a_k=A+kD, ~k=1, 2, 3, \ldots$.
Siegel mentioned explicitly the continued fraction (\ref{natural}), see
\cite[ or p.231 in {\it Gesammelte Abhandlungen} vol. I] {Siegel-1929}.
He obtained  these results as corollaries from the continued fraction expansion of
the ratio of Bessel's functions $J_\lambda(x)$
(see also \cite[formula 9.1.73]{abramowitz+stegun}):
\bee
i\frac{J_{\lambda-1}(2ix)}{J_\lambda(2ix)}=\frac{\lambda}{x}+
\cfrac{1}{\frac{\displaystyle{{\lambda+1}_{\color{white}{9}}}}{\displaystyle{x}^{\color{white}{9} }}+
\cfrac{1}{\frac{\displaystyle{{\lambda+2}_{\color{white}{9}}}}{\displaystyle{x}^{\color{white}{9} }}+
\cfrac{1}{\frac{\displaystyle{{\lambda+3}_{\color{white}{9}}}}{\displaystyle{x}^{\color{white}{9} }}+\ddots}}}
\label{Siegel-ulamek}
\eee
which Siegel has shown to be transcendental for rational $\lambda$ and algebraic
$x\neq 0$ and where the Bessel function of first order is given by
\bee
J_\lambda(x) = \sum_{m=0}^\infty \frac{(-1)^m}{m! \, \Gamma(m+\lambda+1)} {\left({\frac{x}{2}}\right)}^{2m+\lambda}.
\eee
where $\Gamma(x)$ is the Gamma function, see
\cite[first formula on p.231 in {\it Gesammelte Abhandlungen} vol. I] {Siegel-1929}
or  \cite[formula 9.1.10]{abramowitz+stegun}.
For $\lambda=0$ and $x=1$  and taking into account the relation
$J_{-n}(x) = (-1)^n J_{n}(x)$, see e.g. \cite[formula 9.1.5]{abramowitz+stegun},
we get the value of the continued fraction (\ref{natural}):
\bee
s=[0; 1, 2, 3, 4, \ldots]=-i \frac{J_1(2i)}{J_0(2i)}.
\eee

The awkward form (\ref{Siegel-ulamek}) can be written in more pleasant form in terms
of modified Bessel functions of the first kind  $I_\nu(x)$ defined by the series:
\bee
I_\nu(x)=\left(\frac{x}{2}\right)^\nu \sum_{k=0}^\infty \frac{(x/2)^{2k}}{k! \Gamma(\nu+k+1)}.
\eee
There is a following  relation between $I_\nu(x)$ and $ J_\nu(ix)$:
\bee
I_\nu(x)=(-i)^\nu J_\nu(ix),  
\eee
see \cite[formula 9.6.3]{abramowitz+stegun}. Writing $A=\lambda/x, ~~D=1/x$,
 i.e. $\lambda=A/D, ~x=1/D$ we turn (\ref{Siegel-ulamek})
to the more elegant form
\bee
[A; A+D, A+2D, \ldots, A+nD, ...]=\frac{I_{A/D-1}\left(\frac{2}{D}\right)}
{I_{A/D}\left(\frac{2}{D}\right)}
\eee

For $A=-1, ~ D=2$ (or for $\lambda=-1/2$ and $x=1/2$ in (\ref{Siegel-ulamek}))
 we obtain the  value of the continued fraction $s'$ defined by the formula (\ref{odd}):
\bee
s'=[0; 1, 3, 5, 7, 9, \ldots]= 1+\frac{I_{-3/2}(1)}{I_{-1/2}(1)}=
1+i\frac{J_{-3/2}(i)}{J_{-1/2}(i)}.
\label{odd-2}
\eee
The transcendence of (\ref{odd})  follows for $x=i, (i=\sqrt{-1})$ from the
formula known already to Lambert and  Euler \cite{Euler-2005}:
\bee
\tan(x) = \cfrac{x}{1-\cfrac{x^2}{3-\cfrac{x^2}{5-\ddots}}}.
\eee
and the fact that $\tan(x)$ takes transcendental values at algebraic arguments.
From this and from (\ref{odd-2})  as a byproduct we have the identity:
\bee
\tan(i) =i-\frac{J_{-3/2}(i)}{J_{-1/2}(i)}.
\eee

Another  possibility  for partial quotients is the geometrical series : $a_n =q^n$
and we believe that  the corresponding continued fractions:
\bee
G_q=[0; q, q^2, q^3, \ldots]
\eee
are transcendental for all  natural $q \geq 2$.  
This continued fraction is linked to the famous Rogers--Ramanujan
continued fraction  defined by
\bee
RR(q)=\cfrac{q^{1/5}}{1+\cfrac{q}{1+\cfrac{q^2}{1+\cfrac{q^3}{1+\ddots}}}}.
\eee
From the general transformations of continued fractions rules
see \cite[p.9]{Khovanskii} we have the relation:
\bee
q^{4/5}RR(q)=\cfrac{q}{1+\cfrac{q^2}{1+\cfrac{q^3}{1+\cfrac{q^4}{1+\cfrac{q^5}{1+\cfrac{q^6}{\ddots}}}}}}
=\cfrac{1}{q^{-1}+\cfrac{1}{q^{-1}+\cfrac{1}{q^{-2}+\cfrac{1}{q^{-2}+\cfrac{1}{q^{-3}+\cfrac{1}{q^{-3}+\ddots}}}}}}
\eee
In \cite{japonczycy-1996} \cite{japonczycy-1997} it was proved that $RR(q)$
is transcendental for all algebraic $|q|<1$, but it needs some further work to infer
from this the transcendence of $G_q$.
Let us mention, that quite recently K. Dilcher and K. B. Stolarsky \cite{Dilcher2009}
have proved that the continued fractions
\bee
g_q = [q;  q^2, q^4, q^8, \ldots, q^{2^n}, \ldots]
\label{geometr}
\eee
are  transcendental for all integer $q\geq 2$ --- it follows immediately from
(\ref{Davenport-Roth}) and  the double exponential growth:
$Q_n(g_q)>\prod_{k=1}^n q^{2^k}$.     
Adamczewski \cite{Adamczewski2010}  extended this to all complex
$|q|>1$ which are algebraic numbers.
Another (family) class of  transcendental continued fractions  can be found
in \cite{Davison-Shallit-1991}.

Next we can construct a number $f$ where partial quotients are factorials $a_n=n!$:
\bee
f=[0; 1, 2, 6, 24, 120, \ldots, n!, \ldots]=0.6840959001066225003396337\ldots
\label{factorial}
\eee
Even these  partial quotients increase too slowly to apply the Theorem of
Adamczewski and Bugeaud (\ref{Adamczewski-Bugeaud}).  For large $n$ we have
approximately $Q_n(f)\sim \prod_{k=1}^n  k!$. This last product is called superfactorial
and denoted by $n\$$, see also \cite[exercise 4.55]{Knuth1994}.
We prefer the notation $ n!^! = \prod_{k=1}^n  k!$. Superfactorial can be
expressed by the Barnes $G$-function for complex $z$  defined by
\bee
G(z+1)=(2\pi)^{z/2}e^{-(z(z+1)+\gamma z^2)/2}\prod_{n=1}^\infty\left[\left(1+\frac{z}{n}\right)^n e^{-z+z^2/2n}\right].
\eee
It satisfies the functional equation
\bee
G(z+1)=\Gamma(z)G(z)
\eee
and from this we have that
\bee
n!^!=G(n-2).
\eee
The analog of the Stirling formula for $G(z)$ gives \cite{Voros-1987}:
\bee
\log G(z+1)=z^2\left(\frac{\log(z)}{2}-\frac{3}{4}\right)+\frac{z}{2}\log(2\pi)-\ldots
\eee
From this we obtain
\bee
n!^! \sim e^{n^2(\log(n)/2-3/4)}
\eee
and unfortunately
\bee
\frac{\sqrt{\log(Q_n(f))}\log \log (Q_n(f))}{n}\to 0
\eee
hence we do not get transcendentality of $f$ via the Theorem of Adamczewski and Bugeaud.

The continued fraction build  from Fibonacci numbers $a_n=F_n$
\bee
F=[0; 1, 1, 2, 3, 5, 8, \ldots, F_n, \ldots]=0.588873952548933507671231121246787384\ldots  .
\eee
appears at the Sloane The On-Line Encyclopedia of Integer Sequences
as the entry  A073822.  

Apparently both $f$ and $F$ also should be transcendental, but we are not aware
of  the proof of this fact.  The factorial over Fibonacci numbers
behaves as
\bee
\prod_{k=1}^n F_k =\phi^{n(n+1)/2}5^{-n/2}C+\mathcal{O}\left(\phi^{n(n-3)/2}5^{-n/2}\right)
\eee
where $\phi=(1+\sqrt{5})/2$ and $C\approx 1.226742$, see \cite[Exercise 9.41]{Knuth1994}
and it is too slow to use the Davenport -- Roth Theorem.

Let us quote at the end of this Section the following
remarks from the  \cite[p. 104]{Baker-book}: ``And the latter recalls to mind
another outstanding question in Diophantine approximation, namely
whether every continued fraction with unbounded partial quotients is
necessarily transcendental; this too seems very difficult''.
Now there is a common believe that also algebraic numbers of degree $\geq 3$ have
unbounded partial quotients, see e.g.  \cite{Shallit-1992}, \cite{Adamczewski2007}.

\section{Transcendence of $u, u_2, u_d, u_{\FI}, u_q, u_M, u_{r \pm}$}

Because all considered above continued fractions are non-periodic (if there exist
really  infinity of twins, Mersenne primes etc)  they can not
be solutions of a polynomial equations with rational coefficients of the
degree 2, but we believe this statement remains  true for rational
polynomials of  all degrees.
Namely we are convinced  that all  considered above  continued fractions
$u, u_2, u_d, u_{\FI}, u_q, u_M, u_{r \pm}$ are transcendental, however we were not
able to prove it and this problem seems to be extremely difficult.  But if
say $u$ or $u_2$ is not transcendental what the particular polynomial equation
with very special (mysterious) integer coefficients should it satisfy?

It is well known that the Champernowne constant \cite{Champernowne} built by
concatenating consecutive numbers in the base $b$ is transcendental:
\bee
C_b = (\gamma_1)_b(\gamma_2)_b(\gamma_3)_b(\gamma_4)_b\ldots
\eee
where $(\gamma_k)_b$ denotes  number  $k$ expressed in the base $b$ (e.g. in the
common in  computer science notation the twelfth  number in the hexadecimal
base $b=16$ is denoted $C$). In the human base $b=10$ the $C_{10}$  is given by:
\[
C_{10}=0.12345678910111213141516171819202122232425262728293031\ldots
\]
Transcendentality of $C_b$ is the corollary  from the theorem proved by Kurt Mahler
in paper \cite{Mahler-1937} published in 1937. In fact
in this paper \cite{Mahler-1937} Mahler has proved more general result: the
number $\sigma$ obtained by concatenating the values of the positive, integer-valued
increasing polynomial $f(k)$ in the base $b$:
\[
\sigma=f(1)_b f(2)_b f(3)_b \ldots
\]
is transcendental, where $f(k)_b$ denotes the digits of the value of $f(k)$ in
the base system $b$. The case of Champernowne constant is not mentioned in
\cite{Mahler-1937} explicitly but it follows for $f(k)=k$. Let us remark, that
the continued fraction expansion of $C_{10}$ behaves very erratically, with
sporadic partial quotients of  enormous size, for example the 19-th term
is of the order $10^{169}$, 
what is the typical behaviour for the Liouville numbers,
i.e. such numbers that for each $n$ there will be infinity of rationals $A/B$
such that $|C_{10}-A/B|<1/B^n$, see spikes in the Fig. 8. 

The number  $C_{CE}$  obtained by concatenation
of $0.$ with the base 10 representations of the prime numbers in order
\[
C_{CE}=  0.235711131719232931374143475153\ldots
\]
is known as  Copeland–--Erd{\"o}s constant \cite{Copeland-Erdos}. In this paper
Arthur Herbert Copeland and Paul  Erd{\"o}s have shown that $C_{CE}$ is normal,
but apparently it is not proved that  $C_{CE}$ is transcendental.

We have mentioned in the Sect.3  that  $s=[0; 1, 2, 3, 4, 5, 6, \ldots] $
is transcendental. Thus we have the correspondence $s \leftrightarrow C_{b}$
and $u \leftrightarrow C_{CE}$, where both elements of the former pair are shown
to be transcendental and both members of the latter pair are conjectured to be
transcendental.  Of course we have $p_n>n$,  $t_n>n$,  $q_n>n$ etc. but we do
not know how the  transcendence of $u, ~u_2, ~u_q$ follows from these inequalities.

One of the transcendence criterion is the Thue–-Siegel-–Roth Theorem, which we
recall here in the  following form:

{\bf Thue–-Siegel-–Roth Theorem}:  If there exist such $\epsilon>0$ that for
infinitely many fractions $A_n/B_n$   the inequality
\bee
\left|r-\frac{A_n}{B_n}\right|< \frac{1}{B_n^{2+\epsilon}}~~~~~~n = 1,2,3,...,
\label{Thue–-Siegel-–Roth}
\eee
holds, then $r$ is transcendental.

Let us stress, that $\epsilon$ here does not
depend on $n$  ---   it has to be the same for all fractions $A_n/B_n$.
This theorem
suggests the following definition of the measure of irrationality $\mu(r)$:
For a given real number $r$ let us consider the set $\Delta$ of all such exponents
$\delta$ that
\bee
0<\left|r-\frac{P}{Q}\right|<\frac{1}{Q^\delta}
\label{measure-delta}	
\eee
has at most finitely many solutions $(P, Q)$ where $P$ and $Q>0$ are integers.
Then $\mu(r)=\inf_{\delta \in \Delta}\delta$ is called the irrationality measure of $r$
(sometimes any $\delta$ fulfilling (\ref{measure-delta}) is called
irrationality measure and then the smallest $\delta = \mu$ is called the irrationality
exponent).
If the set $\Delta$  is empty, then $\mu(r)$ is defined to be $\infty$ and
$r$ is called a Liouville number. If $r$ is rational then $\mu(r)=1$  and
if $r$ is algebraic of degree $\geq  2$ then $\mu(r)=2$ by the Thue–-Siegel-–Roth Theorem.
There exist real numbers of arbitrary irrationality measure $2\leq \mu < \infty$.  
Namely, the number
\bee
\lfloor a \rfloor +\cfrac{1}{\lfloor a^b \rfloor +\cfrac{1}{\lfloor a^{b^2} \rfloor+
 +\cfrac{1}{\lfloor a^{b^3} \rfloor +\ddots}}}
\eee
where $a>1, ~~b=\mu-1$, has the  irrationality measure $\mu$, see \cite{Brisebarre-2003}.
For the constant $e=\lim_{n\to \infty}(1+1/n)^n$, which has the continued fraction
of a regular form:
\bee
e = [2; 1, 2, 1, 1, 4, 1, 1, 6, 1, 1, 8, 1, 1, 10, 1, 1, 12, 1, 1, \dots, 1,1, 2n, \ldots],
\eee
it is known that $\mu(e)=2$, see \cite[pp.362-365]{Borwein_AGM_PI}.
For $\pi$ it is known that $2\leq \mu(\pi)<7.6304$, see	\cite{Salikhov-2008},
and it is conjectured  \cite[p.203]{Borweiny1989} that $\mu(\pi)=2$.
There is a bound $\delta(n)>2$  for infinitely many $n$ following from the fact
that of any two consecutive convergents to $r$  at least one satisfies
the inequality
\bee
\left|r-\frac{P_n}{Q_n}\right|< \frac{1}{2Q_n^2},
\label{bound-1}
\eee
see \cite[Theorem 18]{Khinchin} or  \cite[Theorem 183]{H-W} and further: of any three
consecutive convergents to $r$,
one at least satisfies
\bee
\left|r-\frac{P_n}{Q_n}\right|< \frac{1}{\sqrt{5} Q_n^2},
\label{bound-2}
\eee
see \cite[Theorem 20]{Khinchin} or \cite[Theorem 195]{H-W}.  Thus writing for convergents satisfying
(\ref{bound-1})  or (\ref{bound-2}) appropriately $\epsilon(n)=\log (2)/\log Q_n$
and  $\epsilon'(n)=\log(5)/2\log Q_n$  the inequality appearing in the
Thue–-Siegel-–Roth Theorem will be satisfied for a given specific $n$.
Of course fractions $P_n/Q_n$ constructed in this way will have
$\lim_{n \to \infty} \epsilon(n)=0$, because $Q_n$ increase monotonically
and there will  be no exponent of $Q_n$ on the r.h.s. of
(\ref{Thue–-Siegel-–Roth}) strictly larger than 2 and  common for all $n$.
In fact, Khinchin \cite{Khinchin} has proved that almost all reals $r$
have $\mu(r)=2$.

The partial quotients of $u, u_d, u_q, u_{\FI}$  grow too slow to use the
Davenport---Roth Theorem,  but if the  behaviour
of the Mersenne primes $\mathcal{M}_n \sim 2^{2^{n e^{-\gamma}}}$ mentioned at the
end of Sect. 6 is valid, then we obtain for large $n$
\bee
Q_n > 2^{c2^{(n+1)e^{-\gamma}}},~~~~~~~~c=\frac{1}{2^{e^{-\gamma}}-1}=2.101893933\ldots
\label{nierownosc}
\eee
and transcendence of $u_\mathcal{M}$ will follow from the Davenport--Roth Theorem
(\ref{Davenport-Roth}). We illustrate the inequality (\ref{nierownosc}) in the Figure 3 --- the values of labels on the $y$--axis give an idea of the order of
$Q_n(\mathcal{M}_n)$: the largest for $n=47$ is of the order
$Q_{47}= e^{1.9984\ldots \times 10^8}=2.32928\ldots \times 10^{86789810}$!

\begin{figure}
\vspace{-2.3cm}
\hspace{-3.5cm}
\begin{center}
\includegraphics[height=0.4\textheight,angle=0]{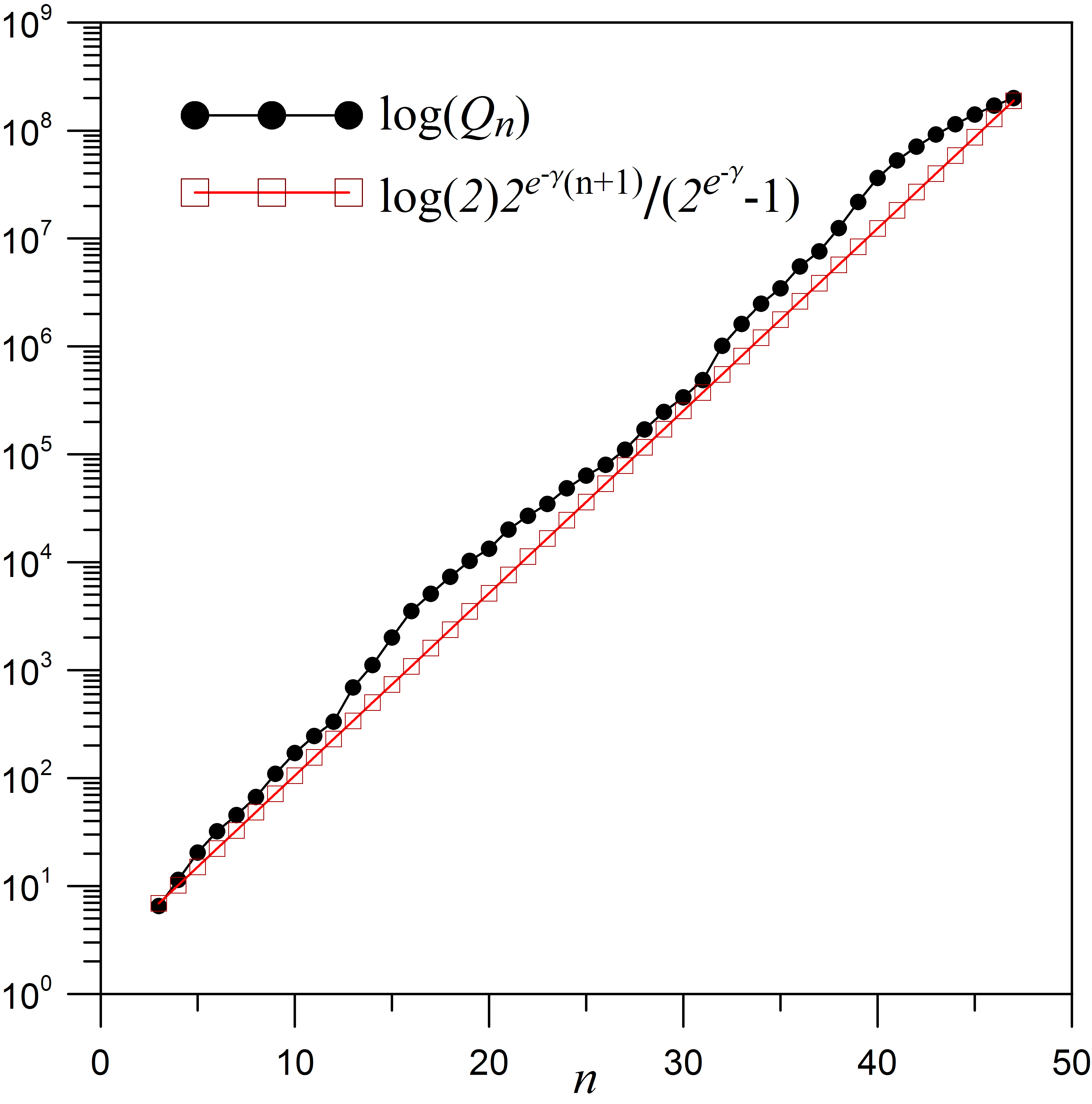} \\
\vspace{0.2cm} {\small Fig.3 Illustration of the inequality (\ref{nierownosc})
for $3\leq n \leq 47$. Although the last points seem to coincide in fact
$Q_{47}= 2.32928\ldots \times 10^{86789810}$, while $2^{c2^{48e^{-\gamma}}}=
1.21513\ldots \times 10^{82034318}$ --- hundreds thousands orders of difference!}
\end{center}
\end{figure}

Usually the number $r$ in question (for example $e, \pi, \zeta(3)$, etc.)
is given by some definition not involving continued fractions, but here we have
expressions of $u, \ldots u_{\mathcal{M}}$ {\bf only} by continued fractions and
we can not calculate directly the differences $|r-{P_n}/{Q_n}|$, like it is possible
for example for Liouville transcendental  numbers or for $e$. For this last case,
as mentioned earlier,
the possibility of explicit calculation of the difference $|e-{P_n}/{Q_n}|$
gives that $\mu(e)=2$ see \cite[pp. 351--371]{Borwein_AGM_PI}.
 We do not have any ideas now how
to express $g_A, f, u, ~u_d, ~u_q, ...$  independently by means  of formulas not
involving continued fractions.
Nevertheless we have made the plot of the exponent $\delta(n)$ in the difference:
\bee
\left|U-\frac{P_n}{Q_n}\right|=  \frac{1}{Q_n^{\delta(n)}}~~~~~~n = 1,2,3,...,
\eee
where $U$ stands for $u, ~u_d, ~u_q, ...$ and $P_n/Q_n$ are convergents of
continued fractions for $U$ --- it is a well known fact that convergents of
continued fractions are the best rational approximations.

In the Figures 4--8 we present plots
of $\delta(n)$ for $u, u_q$ and for $\pi$, $e$ as well as  for $C_{10}$
for comparison.  First we have calculated
$u, ~u_q, ...$  with  150000 digits accuracy from the generic definition
by constructing the continued fractions with a many thousands partial denominators.
Next we have calculated $P_n$ and $Q_n$ for $n$ until the difference $|U-P_n/Q_n|$ was
zero in prescribed accuracy. From the differences $|U-\frac{P_n}{Q_n}|$
we calculated the values of $\delta(n)=-\log|U-\frac{P_n}{Q_n}|/\log(Q_n)$ and
the sample of results is plotted on Fig.4 and 5 for $u$ and $u_q$.
 The bound following from (\ref{bound-1})
is fulfilled for all $n<2100$. In the next Figures we present the plot
of $\delta(n)$  for $\pi$ (Fig.6), $e$  (Fig.7)  and $C_{10}$ (Fig.8). For $C_{10}$ we
have plotted
$\delta(n)-2$ because values of this difference changes by many orders, in contrast to
smooth behavior seen in the Figs. 4-7. The spikes seen in the  Fig.8 are
similar to the behavior of the Liouville transcendental numbers,
but the last statement in  \cite{Mahler-1937} asserts that  $C_{10}$ is not
the Liouville number.

In \cite{Sondow-2004} J. Sondow has proved that:
\bee
\mu(r)=1+\limsup_{n \to \infty}\frac{\log Q_{n+1}}{\log Q_n}=2+\limsup_{n \to \infty}\frac{\log a_{n+1}}{\log Q_n}
\eee

From this we have for $u$ as $a_n=p_n\sim n\log(n)$ and  for large $n$
$Q_n \sim n^n$ that $\mu(u)=2$
and the same for $u_2, u_q$, but for Mersenne primes we get from the Wagstaff
conjecture:
\bee
\mu(u_\mathcal{M})<2+2^{e^{-\gamma}}-1=2.47477\ldots
\eee
But if there is only finite number of Mersenne primes (and hence finitely
many even perfect numbers), then  $\mu(u_\mathcal{M})=1$. In the Fig.9 we present the
plot of $\delta(\mathcal{M}_n; n)=-\log|u_\mathcal{M}-P_n/Q_n|/\log(Q_n)$ and indeed
the values oscillate around $1+2^{e^{-\gamma}}=2.47477\ldots$.

We propose the conjecture which indicates the way to deduce the transcendence
of some continued fractions from transcendence of another ones:

{\bf Conjecture $(\star)$}: Let $\alpha=[a_0; a_1, a_2, \ldots], ~~\lim_{n\to \infty} a_n = \infty$,
and $\beta=[b_0; b_1, b_2, \ldots]$,  where  $a_n, b_n \in \mathbb{N}$.  Suppose there
exists such  $n_0$ that for all $n>n_0$ the inequality $b_n>a_n$  holds. If $\alpha$
is  transcendental then $\beta$ is also transcendental.

The condition $ \lim_{n\to \infty} a_n = \infty$ is necessary: if $a_n$ is bounded,
say $a_n < A$ for $\forall n, ~~A\in \mathbb{N}$, then
$\beta=[A; 2A, 2A, \ldots]=\sqrt{1+A^2}$.
Also for transcendental $b_n$ the above conjecture probably is not true. When the
Conjecture  $(\star)$ will be proved it will suffice for our purposes to  invoke the
transcendence of the continued fraction
$s = [0; 1, 2, 3, 4, \ldots]$ (\ref{natural}), as for all examples from Sect.2
we have $a_n>n$ (then also some examples from Sect.3 will be transcendental, like
$f$ with $a_n=n!$ and $F$ with $a_n=F_n$).

\section{Final remarks}

We have raised above some questions concerning the transcendence of 
continued fractions with partial quotients given by prime numbers of a few special
forms.  We hope that the experimental results  reported  above will stimulate
further research in the field.

{\bf Acknowledgement } I would like to thank professors   Boris Adamczewski,
Jaroslav Han\v{c}l and Michael Waldshmidt for e-mail exchange.

\begin{figure}
\vspace{-2.3cm}
\hspace{-3.5cm}
\begin{center}
\includegraphics[height=0.4\textheight,angle=0]{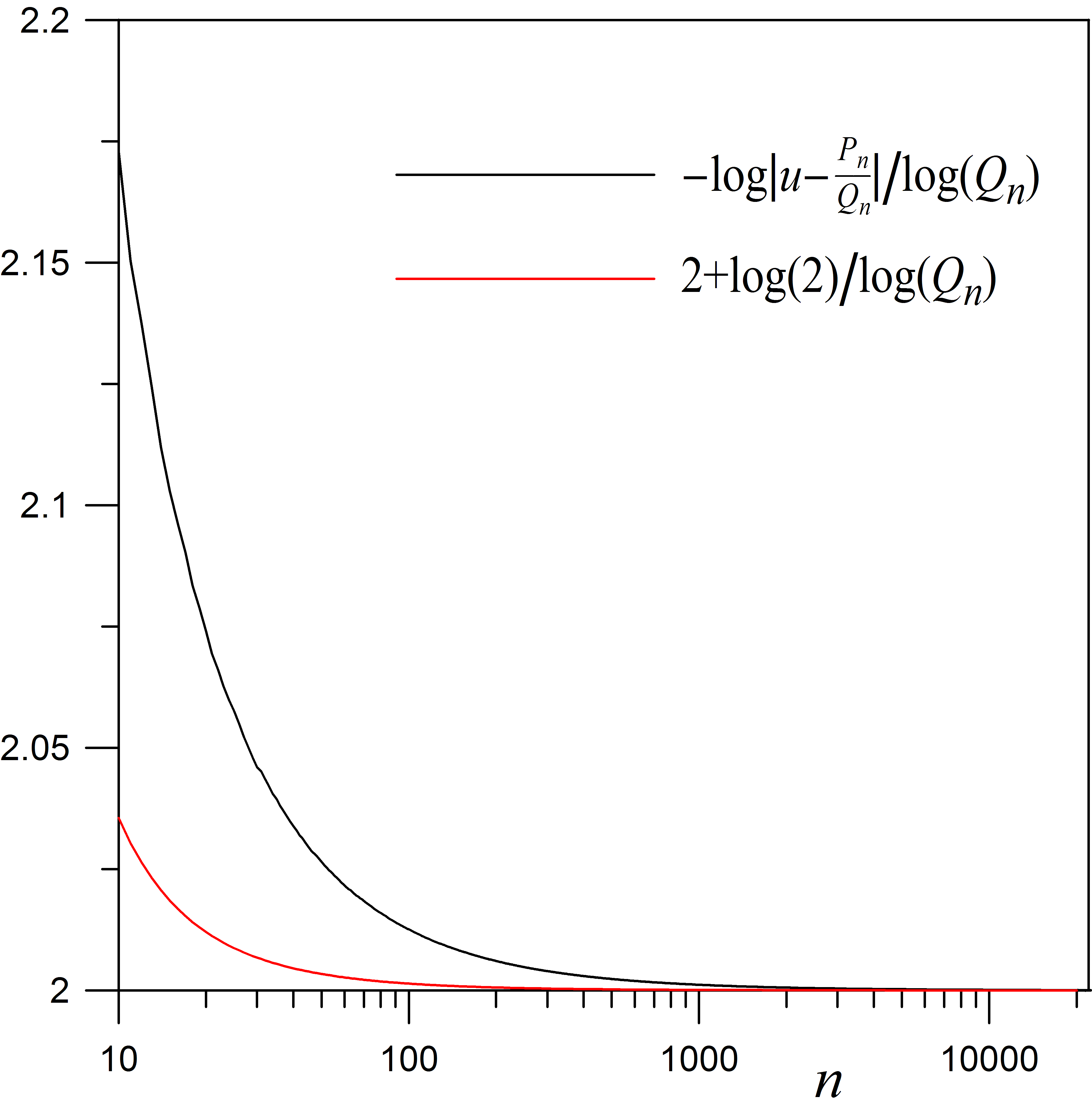} \\
\end{center}
\vspace{0.0cm} {\small  Fig.4 The plot of $\delta(n)$ (black) and the bound
$2+\log (2)/\log Q_n$  (red)  for $u$ following from the (\ref{bound-1})
up to $n=22380$  (computatuions were done in precision 150000 digits
and the value of $|u-P_{22380}/Q_{22380}|$ was zero with accuracy 150000 digits).
For each $n$ the bound (\ref{bound-1}) (as well as bound (\ref{bound-2}))}
is fulfilled.
\end{figure}

\begin{figure}
\vspace{-0.7cm}
\hspace{-3.5cm}
\begin{center}
\includegraphics[height=0.4\textheight,angle=0]{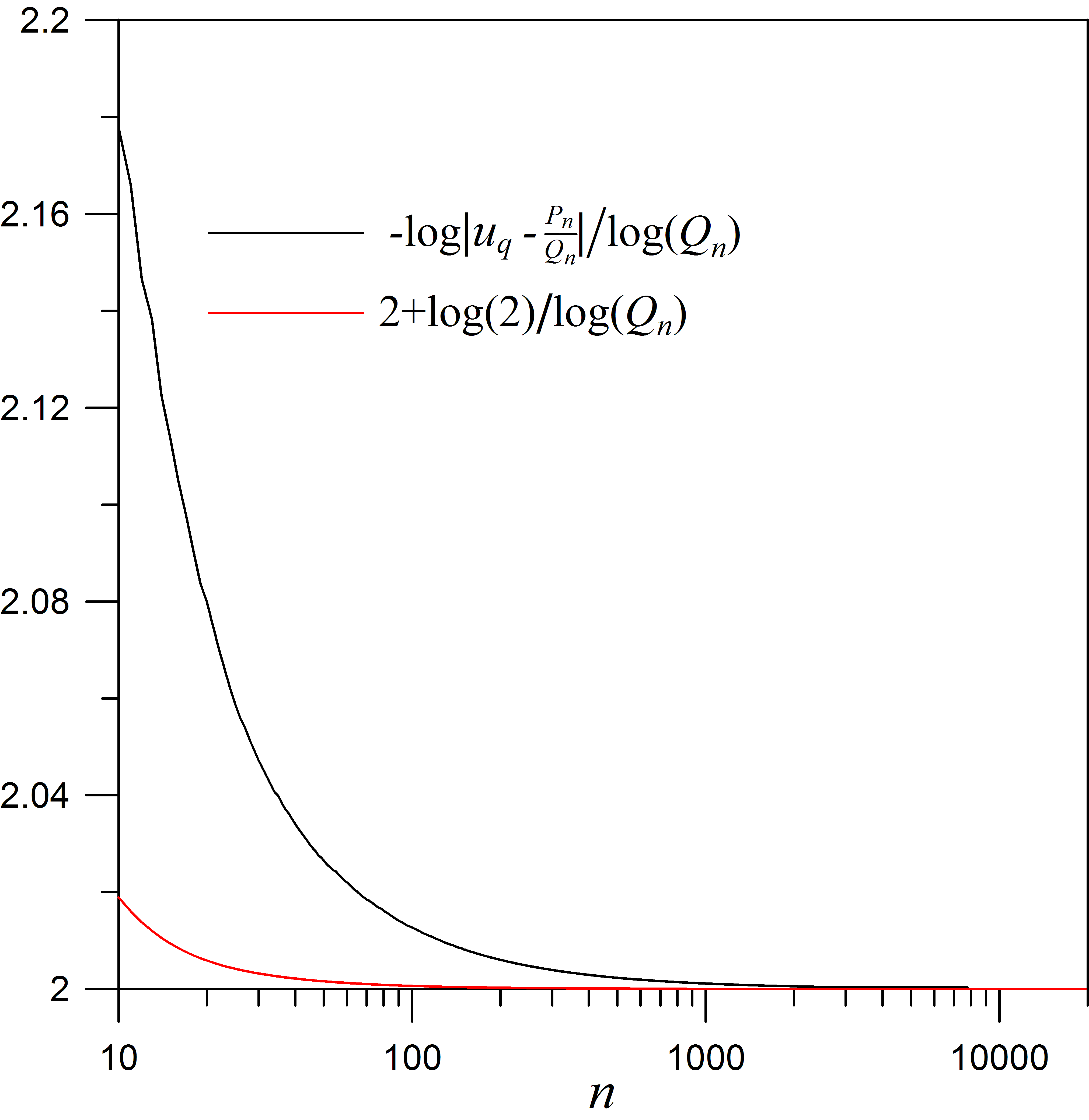} \\
\vspace{0.0cm} {\small Fig.5 The plot of $\delta(n)$ (black) and the bound
$2+\log (2)/\log Q_n$   for $u_q$ following from the (\ref{bound-1}).
For each $n$ the bound (\ref{bound-1}) is fulfilled.}
\end{center}
\end{figure}

\newpage

\begin{figure}[pht]
\vspace{-1cm}
\begin{minipage}{13.8cm}
\hspace{3.5cm}
\begin{center}
\includegraphics[height=0.34\textheight,angle=0]{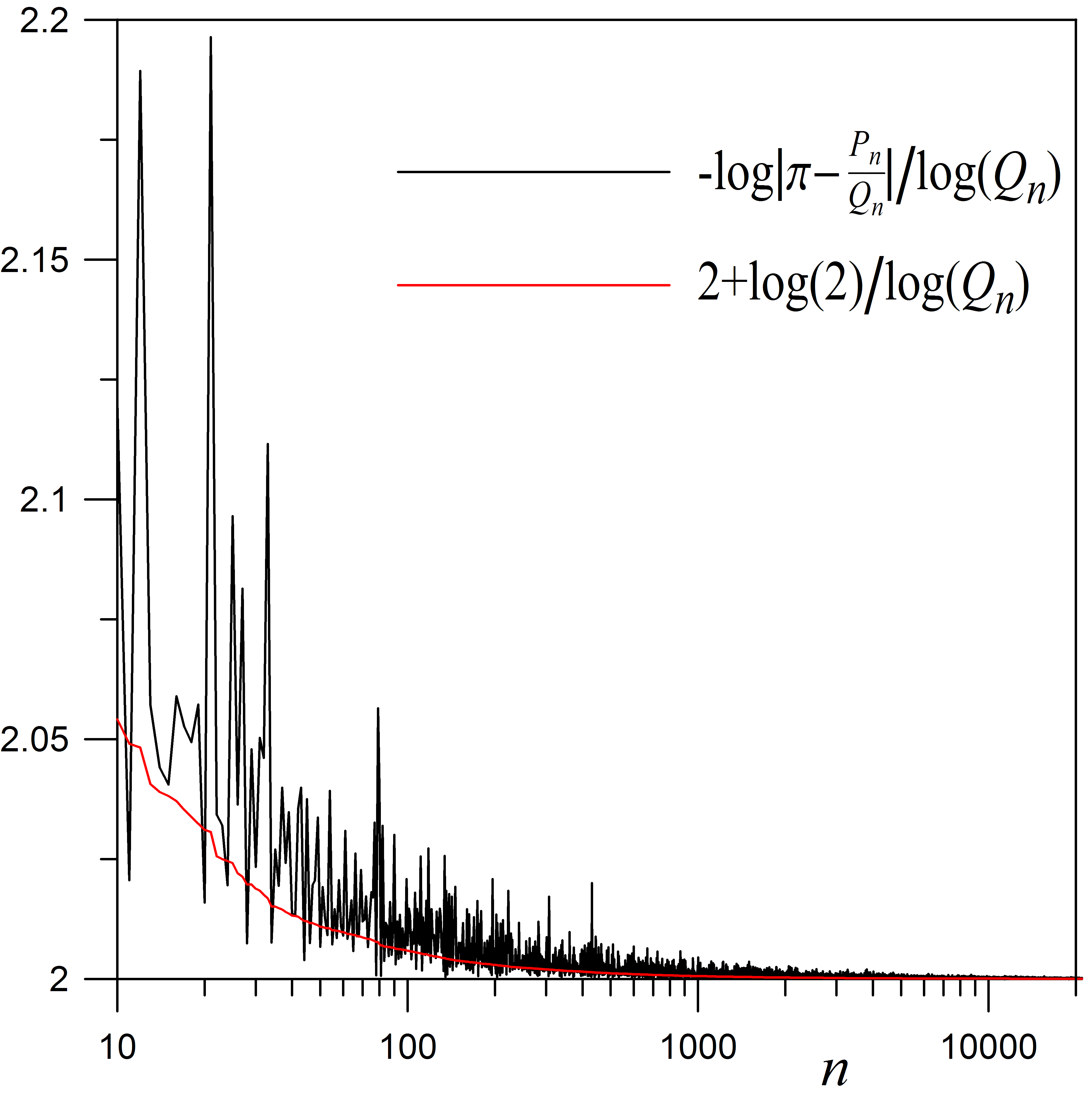} \\
\end{center}
\vspace{0.0cm}{\small  Fig.6 The plot of $\delta(n)$ (black) and the bound
$2+\log (2)/\log Q_n$ (red)   for $\pi $ following from the (\ref{bound-1})
up to $n=20997$  (computatuions were done in precision 50000 digits
and the value of $|\pi-P_{20997}/Q_{20997}|$ was zero with prescribed accuracy).}
\end{minipage}

\vspace{0cm}
\begin{minipage}{13.8cm}
\hspace{3.5cm}
\begin{center}
\includegraphics[height=0.35\textheight,angle=0]{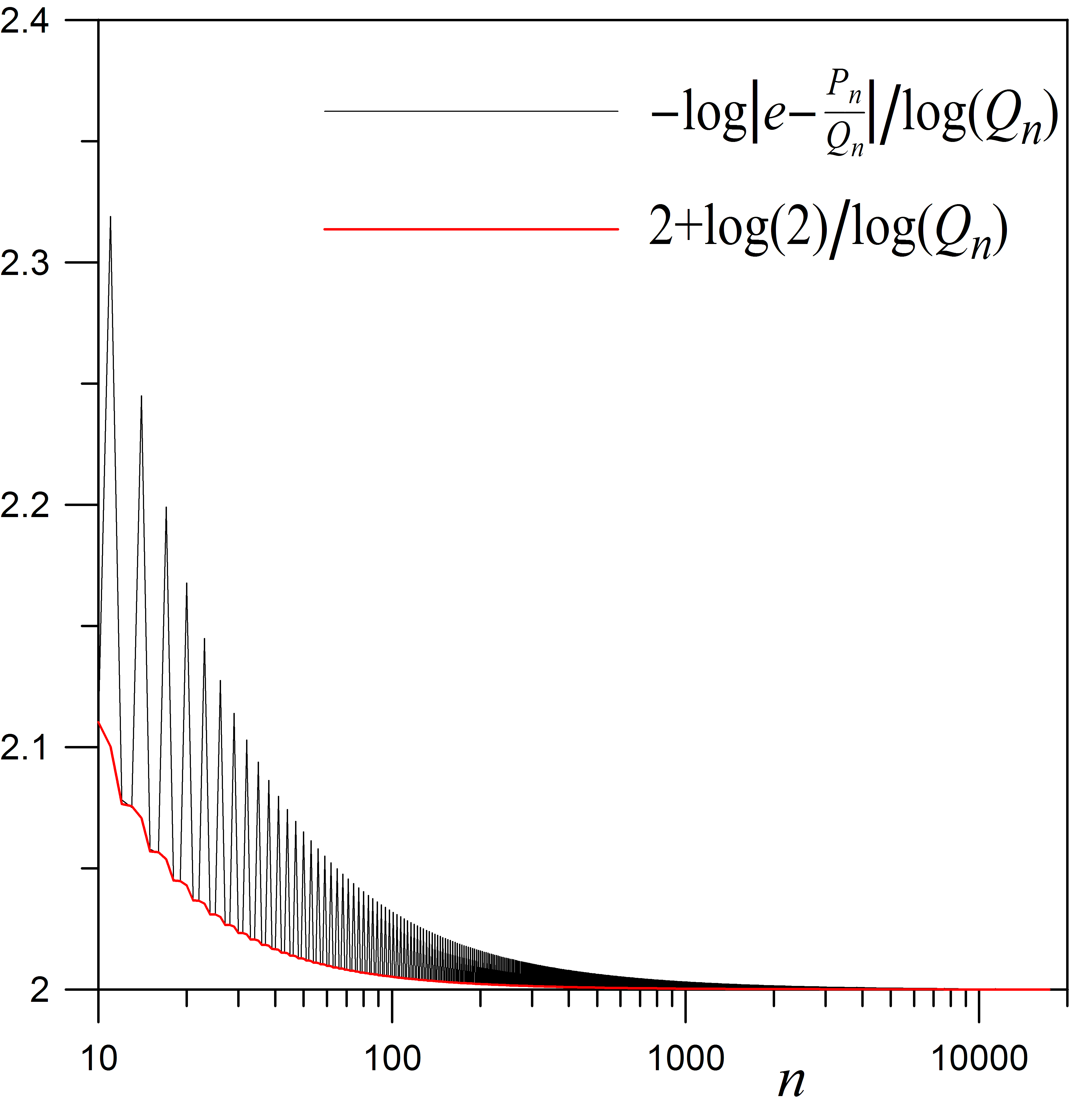} \\
\end{center}
\vspace{0.0cm} {\small Fig.7 The plot of $\delta(n)$ (black) and the bound
$2+\log (2)/\log Q_n$ (red)  for $e$  following from the (\ref{bound-1})
up to $n=17365$ (computatuions were done in precision 45500 digits
and at this $n$ the value of $|e-P_{17365}/Q_{17365}|$ was zero).
The periodic structure of the continued fraction expansion for $e$
is  clearly seen.}
\end{minipage}
\end{figure}

\newpage

\begin{figure}
\vspace{-2.3cm}
\hspace{-3.5cm}
\begin{center}
\includegraphics[height=0.6\textheight, width=0.9\textwidth, angle=0]{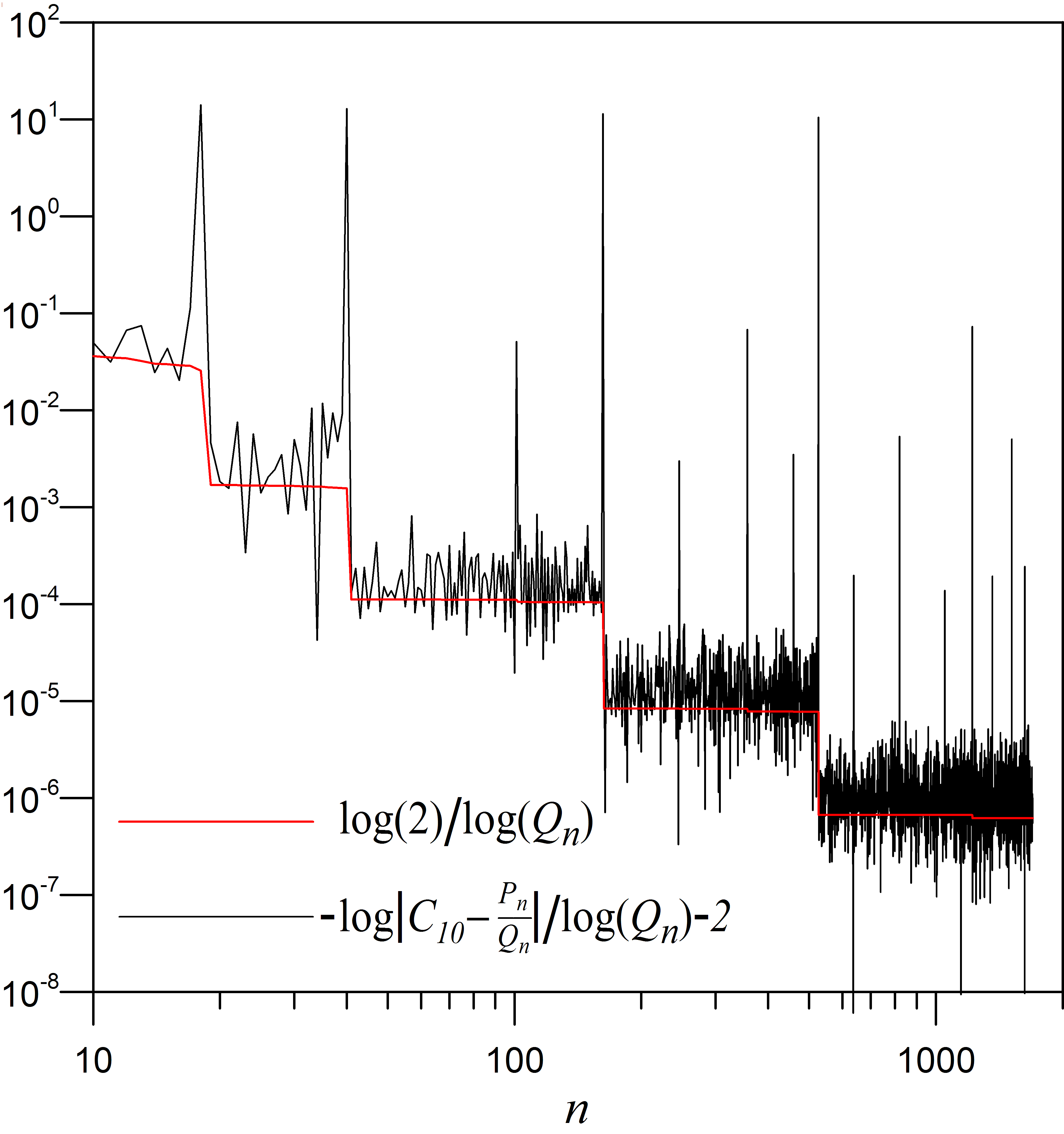} \\
\end{center}
\vspace{-0.4cm}  Fig.8 The plot of $\delta_{C_{10}}(n)-2$ (black) and the bound
$\log (2)/\log Q_n$ (red)  for $C_{10}$ following from the (\ref{bound-1}).
Because of the weird behavior of the partial quotients of continued fraction expansion
for  $C_{10}$ we have subtracted 2 from   $\delta(n)$  and plotted the graph
with the $y$ axis in the logarithmic scale. After each extremely large partial
quotient $a_n$ there is an abrupt drop in the values of $\delta(n)$ and the bound
$ \log (2)/\log Q_n$ with accompanying spike for $n-1$, see (\ref{error}).
It took almost 4 days CPU time to get data for this plot.
Collecting data was done in a few separate runs with different precisions.
Because the partial quotient $a_{526}>10^{411100}$ and $Q_{527}>10^{449994}$
the calculations for $526\leq n < 1708$
was performed with 1,000,000 digits precision, see eq.(\ref{error}).
We stopped at $n=1707$ because $a_{1708}>10^{4911098}$. Spikes of $\delta_{C_{10}}(n)$
many orders higher then neighboring values suggest that $C_{10}$ may be the
transcendental number of Liouville type, but it in contradiction with the last statement
of the paper \cite{Mahler-1937}.
\end{figure}

\newpage

\begin{figure}
\vspace{-2.3cm}
\hspace{-3.5cm}
\begin{center}
\includegraphics[height=0.4\textheight,angle=0]{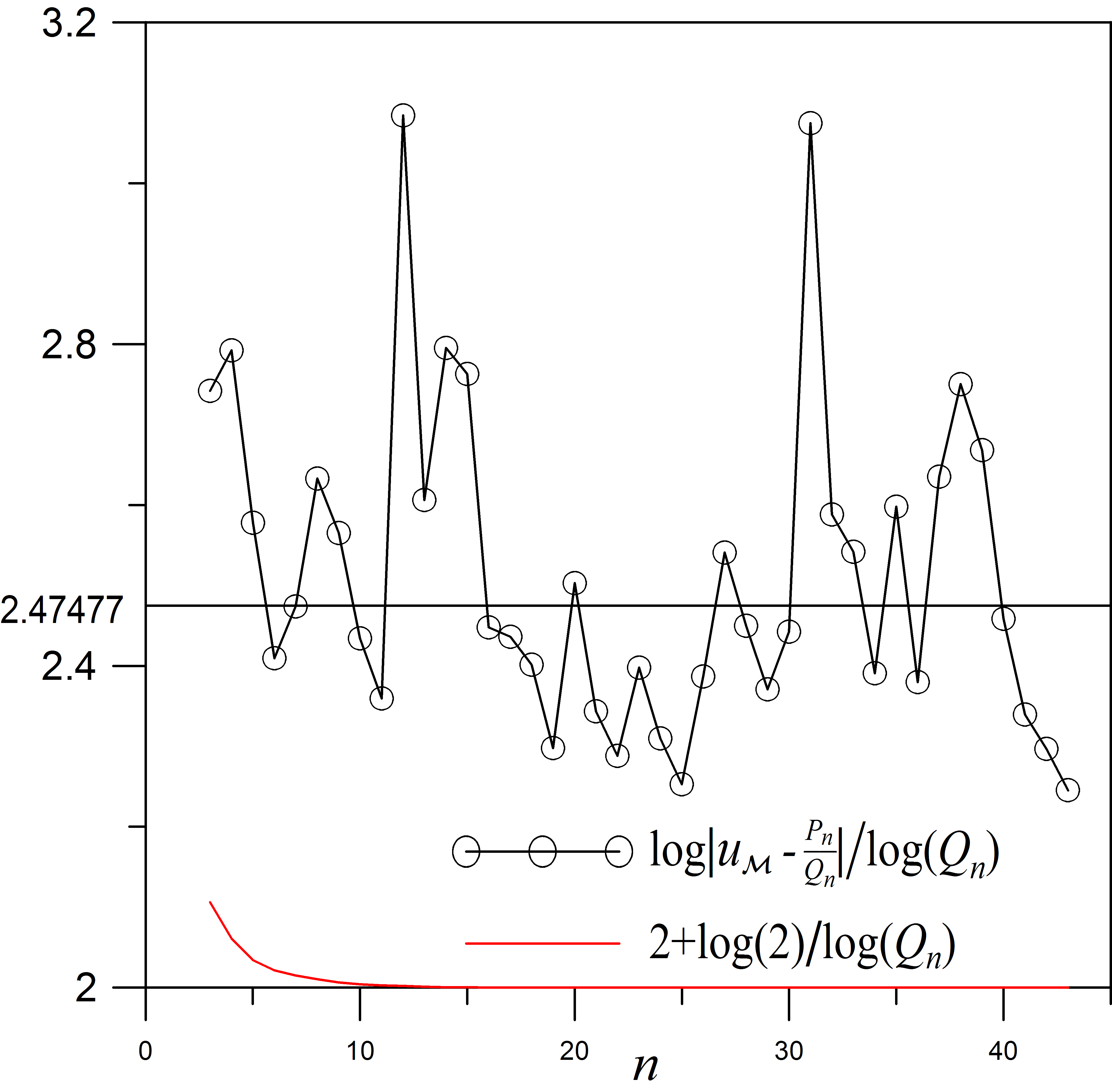} \\
\vspace{0.2cm} {\small Fig.9 The plot of $-\log|u_\mathcal{M}-P_n/Q_n|/\log(Q_n)$ (black)
and the bound $2+\log (2)/\log Q_n$  (red) following from the (\ref{bound-1}) for
$ ~3\leq n \leq 43$.  Here the value of  $u_\mathcal{M}$ was obtained from
all 47 known Mersenne primes with more than 120 millions digits:
the accuracy was better than $10^{-121949117}$.
The denominators $Q_n$ grow very fast and
the bound $2+\log(2)/\log(Q_n)$ tends quickly to 2. It took 12 days CPU time on
the AMD Opteron 2700 MHz processor to collect data for $n\leq 40$: the  point $n=40$
needed precision of almost 40,000,000 digits, as $|u_\mathcal{M}-P_{40}/Q_{40}|
=1.5033\times 10^{-38789567}$, while $1/Q_{40}^2=4.501\ldots \times 10^{-31553835}$}.
To calculate the difference $|u_\mathcal{M}-P_n/Q_n|$ for $n=41, 42, 43$ the precision
of  100000000 digits was needed and one point took 6 days on the same processor,
as for example  $|u_\mathcal{M}-P_{43}/Q_{43}|< 10^{-89770217}$.
\end{center}
\end{figure}


\begin{thebibliography}{10}
\bibitem{Euler-2005}
L.~Euler.
\newblock A commentary on the continued fraction by which the illustrious {L}a
  {G}range has expressed the binomial powers, 2005.
\newblock arXiv.org:math.HO/0507459.
\end{thebibliography}
\end{document}